\documentclass[letterpaper,twoside,final,notitlepage,onecolumn,thmsb]{article}
\usepackage{amssymb}
\begin{document}

\author{David Carf\`{i}}
\title{$^{\mathcal{S}}$Bases in $^{\mathcal{S}}$Linear Algebra}
\date{}
\maketitle

\begin{abstract}
In this paper we define the $^{\mathcal{S}}$bases for the spaces of tempered
distributions. These new bases are the analogous of Hilbert bases of
separable Hilbert spaces for the continuous case (they are indexed by $m$-dimensional Euclidean spaces) and enjoy properties
similar to those shown by algebraic bases in the finite dimensional case. The $^{
\mathcal{S}}$bases are one possible rigorous and extremely manageable mathematical model for the ``physical'' bases used in Quantum Mechanics.
\end{abstract}

\bigskip

\bigskip

\bigskip

\bigskip 

\bigskip 

\section{$^{\mathcal{S}}$\textbf{Linear independence}}

\bigskip

A finite family $v$ of vectors of a vector space is said linearly
independent if any zero linear combination of $v$ is given by a zero system
of coefficients. This is exactly the definition we generalize in the
following.

\bigskip

\textbf{Definition (of }$^{\mathcal{S}}$\textbf{linear independence).}\emph{%
\ Let\ }$v\in \mathcal{S}(\Bbb{R}^{m},\mathcal{S}_{n}^{\prime })$\emph{\ be
an }$^{\mathcal{S}}$\emph{family of tempered distributions in the space }$%
\mathcal{S}_{n}^{\prime }$\emph{\ indexed by the Euclidean space }$\Bbb{R}%
^{m}$\emph{. The family }$v$\emph{\ is said to be }$^{\mathcal{S}}$\emph{%
\textbf{linearly independent}, if any coefficient distribution }$a\in 
\mathcal{S}_{m}^{\prime }$\emph{\ such that} 
\[
\int_{\Bbb{R}^{m}}av=0_{\mathcal{S}_{n}^{\prime }} 
\]
\emph{must be the zero distribution }$0_{\mathcal{S}_{m}^{\prime }}$\emph{.
In other terms the family }$v$\emph{\ is said }$^{\mathcal{S}}$\emph{%
linearly independent if and only if any zero }$^{\mathcal{S}}$\emph{linear
combination of the family }$v$\emph{\ has necessarily a zero coefficient
system.}

\bigskip

We shall denote the superposition of a family of distribution $v$ with
respect to a distributional system of coefficients $a$, also by the
multiplication symbol $a.v$.

\bigskip

\textbf{Example (The Dirac family). }\emph{The Dirac family\ in }$\mathcal{S}%
_{n}^{\prime }$\emph{\ is }$^{\mathcal{S}}$\emph{linearly independent.} In
fact, we have 
\[
\int_{\Bbb{R}^{n}}u\delta =u, 
\]
for all $u\in \mathcal{S}_{n}^{\prime }$, and then the relation $u.\delta
=0_{\mathcal{S}_{n}^{\prime }}$ is equivalent to $u=0_{\mathcal{S}%
_{n}^{\prime }}$.

\bigskip

\textbf{Example (the Fourier families).}\ Let $a,b$\ be two real non-zero
numbers, recall that the $(a,b)$-Fourier\ family in the space of tempered
distributions $\mathcal{S}^{\prime }(\Bbb{R}^{n},\Bbb{C})$ is the family 
\[
\left( a^{-n}\left[ e^{-ib(p\mid \cdot )}\right] \right) _{p\in \Bbb{R}%
^{n}}, 
\]
of smooth and bounded regular tempered distributions (recall that $[f]$
denote the distribution canonically associated to a locally summable
function $f$) . In the particular case $a=1$\ and $b=-1/\hbar $\ (with $%
\hbar $ the reduced Planck constant) we obtain what we call the De Broglie
family, i.e. the family 
\[
\left( \left[ e^{(i/\hbar )(p\mid \cdot )}\right] \right) _{p\in \Bbb{R}%
^{n}}. 
\]
\emph{All the Fourier families are }$^{\mathcal{S}}$\emph{linearly
independent. }In fact, let $\varphi $ be the $\left( a,b\right) $-Fourier
family, and let 
\[
\int_{\Bbb{R}^{n}}u\varphi =0_{\mathcal{S}_{n}^{\prime }(\Bbb{C})}. 
\]
For every test function $\phi $ in $\mathcal{S}_{n}(\Bbb{C})$, we have 
\begin{eqnarray*}
0 &=&\left( \int_{\Bbb{R}^{n}}u\varphi \right) \left( \phi \right) = \\
&=&u\left( \widehat{\varphi }\left( \phi \right) \right) = \\
&=&u\left( \mathcal{S}_{(a,b)}\left( \phi \right) \right) = \\
&=&\mathcal{F}_{(a,b)}\left( u\right) \left( \phi \right) ,
\end{eqnarray*}
where $\mathcal{S}_{(a,b)}$ is the Fourier-Schwartz transformation on the
space $\mathcal{S}_{n}(\Bbb{C})$, so that the Fourier transform of the
distribution $u$ must be zero, i.e. 
\[
\mathcal{F}_{(a,b)}\left( u\right) =0_{\mathcal{S}_{n}^{\prime }(\Bbb{C})}, 
\]
and hence $u=0_{\mathcal{S}_{n}^{\prime }(\Bbb{C})}$, being the Fourier
transformation $\mathcal{F}_{(a,b)}$ an injective operator.

$\bigskip $

\section{\textbf{Linear and }$^{\mathcal{S}}$\textbf{linear independence}}

\bigskip

A first elementary connection between the classic linear independence and
the new $^{\mathcal{S}}$linear independence is given by the following
theorem, which states that the $^{\mathcal{S}}$linear independence is a
stronger requirement for an $^{\mathcal{S}}$family than the requirement of
the simple linear independence.

\bigskip

Recall that an infinite family $v$ of vectors of a vector space is said to
be linearly independent if any finite subfamily of $v$ is linearly
independent.

\bigskip

\textbf{Theorem.}\emph{\ Let }$v\in \mathcal{S}(\Bbb{R}^{m},\mathcal{S}%
_{n}^{\prime })$\emph{\ be an }$^{\mathcal{S}}$\emph{linearly independent
family. Then, the family }$v$\emph{\ is linearly independent. Consequently,
the }$^{\mathcal{S}}$\emph{linear hull }$^{\mathcal{S}}\mathrm{span}\left(
v\right) $\emph{\ is an infinite dimensional subspace of }$\mathcal{S}%
_{n}^{\prime }$\emph{.}

\bigskip

\emph{Proof.} By contradiction, assume the family $v$ linearly dependent.
Then there exists a linearly dependent finite subfamily of $v$. More
precisely, there are a positive integer $k\in \Bbb{N}$ and a $k$-sequence $%
\alpha \in \left( \Bbb{R}^{m}\right) ^{k}$ of points belonging to the $m$%
-dimensional Euclidean space $\Bbb{R}^{m}$ such that the $k$-sequence of
distributions $v_{\alpha }=(v_{\alpha _{i}})_{i=1}^{k}$, extracted from the
family $v$ by means of the index selection $\alpha $ is a linearly dependent
system of $\mathcal{S}_{n}^{\prime }$. Then there exists a non-zero scalar $%
k $-tuple $a\in \Bbb{K}^{k}$ such that the finite linear combination $\sum
av_{\alpha }$ is zero, that is such that 
\[
\sum_{i=1}^{k}a_{i}v_{\alpha _{i}}=0_{\mathcal{S}_{n}^{\prime }}. 
\]
Consider the tempered distribution $d=\sum_{i=1}^{k}a_{i}\delta _{\alpha
_{i}}$ as a coefficient distribution, we have 
\begin{eqnarray*}
\int_{\Bbb{R}^{m}}dv &=&\int_{\Bbb{R}^{m}}\left( \sum_{i=1}^{k}a_{i}\delta
_{\alpha _{i}}\right) v= \\
&=&\sum_{i=1}^{k}a_{i}\left( \int_{\Bbb{R}^{m}}\delta _{\alpha _{i}}v\right)
= \\
&=&\sum_{i=1}^{k}a_{i}v_{\alpha _{i}}= \\
&=&0_{\mathcal{S}_{n}^{\prime }}.
\end{eqnarray*}
Now, since the distribution $d$ is different from the zero distribution $0_{%
\mathcal{S}_{m}^{\prime }}$, the preceding equality contradicts the $^{%
\mathcal{S}}$linear independence of the family $v$, against our assumptions. 
$\blacksquare $

$\bigskip $

\section{\textbf{Topology and }$^{\mathcal{S}}$\textbf{linear independence}}

\bigskip

The last theorem of the above section shows that, for what concerns the $^{%
\mathcal{S}}$families, the $^{\mathcal{S}}$linear independence implies the
usual linear independence. Actually, the $^{\mathcal{S}}$linear independence
is more restrictive than the linear independence, as we shall see later by a
simple notable example of family which is linearly independent but not $^{%
\mathcal{S}}$linearly independent. On the contrary the $^{\mathcal{S}}$%
linear independence is less restrictive than the $\beta (\mathcal{S}%
_{n}^{\prime })$-topological independence, as it is shown below.

\bigskip

\textbf{Topological independence.} We recall that a system of vectors $%
v=(v_{i})_{i\in I}$ in the space $\mathcal{S}_{n}^{\prime }$, indexed by a
non-void index set $I$, is said $\beta (\mathcal{S}_{n}^{\prime })$-\emph{%
topologically free} (respectively, $\sigma (\mathcal{S}_{n}^{\prime })$%
-topologically free) if and only if there exists a family $L=(L_{i})_{i\in
I} $ of $\beta (\mathcal{S}_{n}^{\prime })$-continuous (respectively, $%
\sigma (\mathcal{S}_{n}^{\prime })$-continuous) linear forms on $\mathcal{S}%
_{n}^{\prime }$ such that 
\[
L_{i}(v_{k})=\delta _{ik}, 
\]
for any pair $(i,k)\in I^{2}$, where the family $\delta =(\delta
_{ik})_{(i,k)\in I^{2}}$ is the Kronecker family on the square $I^{2}$. Note
that the above relation can be written as $L\otimes v=\delta $, where $%
\delta $ is the Kronecker family.

If the family $v$ is not topologically free it is said \emph{topologically
bound}. If the family $v$ is topologically free, any family $L$ of
continuous linear forms, satisfying the above relations is said a \emph{dual
family of the family} $v$. So, to say that a family $v$ is topologically
free is equivalent to say that $v$ has a dual family of linear continuous
forms.

Recalling that (by reflexivity) any continuous linear functional on the
space $\mathcal{S}_{n}^{\prime }$ is canonically and univocally
representable by a test function in $\mathcal{S}_{n}$, to say that the
family $v$ is topologically free is equivalent to say that the
bi-orthonormality condition 
\[
\langle g_{i},v_{k}\rangle =\delta _{ik}, 
\]
is true, for any pair $(i,k)\in I^{2}$, where $g$ is a suitable family of
test functions.

\bigskip

\textbf{Theorem.} \emph{Every }$^{\mathcal{S}}$\emph{family in the space }$%
\mathcal{S}_{n}^{\prime }$\emph{\ is }$\beta (\mathcal{S}_{n}^{\prime })$%
\emph{-topologically bound and, thus, }$\sigma (\mathcal{S}_{n}^{\prime })$%
\emph{-topologically bound. Consequently, no }$^{\mathcal{S}}$\emph{family
has a dual family of test functions.}

\bigskip

\emph{Proof.} Let $v$ be any $^{\mathcal{S}}$family in the space $\mathcal{S}%
_{n}^{\prime }$ indexed by $\Bbb{R}^{m}$. And let $L$ be an arbitrary family
in the dual $\mathcal{S}_{n}^{\prime \prime }$ indexed by the same index
set. Being the Schwartz space $\left( \mathcal{S}_{n}\right) $ reflexive,
for every $i$, there is a test function $g_{i}$ in $\mathcal{S}_{n}$
canonically generating the functional $L_{i}$, that is such that 
\[
L_{i}=\left\langle .,g_{i}\right\rangle . 
\]
In other terms, the test function $g_{i}$ is such that $L_{i}(u)=u(g_{i})$,
for every tempered distribution $u$ in $\mathcal{S}_{n}^{\prime }$. Assume
the existence of an index $i$ such that $L_{i}(v_{i})=1$, then we deduce 
\[
1=L_{i}(v_{i})=v_{i}(g_{i})=v(g_{i})\left( i\right) , 
\]
being $v$ an $^{\mathcal{S}}$family, the function $v(g_{i})$ is continuous,
then there is a neighborhood $U$ of the point $i$ in which the function $%
v(g_{i})$ is strictly positive. Then, for every point $k$ in the
neighborhood $U$, we have 
\[
L_{i}(v_{k})=v_{k}(g_{i})=v(g_{i})\left( k\right) >0, 
\]
and then $L$ cannot verify the condition 
\[
L_{i}(v_{k})=\delta _{ik}, 
\]
for any pair $(i,k)\in I^{2}$. So we cannot find a dual family of
functionals for $v$ and, consequently, $v$ cannot be topologically
independent. $\blacksquare $

\bigskip

\textbf{Note.} By the same proof, it is possible to prove that \emph{every }$%
C^{0}$\emph{-family of tempered distributions is strongly topologically bound%
}. Consequently every smooth family of tempered distributions is strongly
topologically bound.

$\bigskip $

\section{\textbf{Multiplicity of representations}}

\bigskip

It's simple to prove the following property that characterizes the $^{%
\mathcal{S}}$linear dependence of a family of distributions by explicit
multeplicity of representations of some member of the family itself.

\bigskip

\textbf{Property.}\emph{\ An }$^{\mathcal{S}}$\emph{family }$v$\emph{\ in }$%
\mathcal{S}_{n}^{\prime }$\emph{\ indexed by }$\Bbb{R}^{m}$\emph{\ is }$^{%
\mathcal{S}}$\emph{linearly dependent if and only if there is a point index }%
$p$\emph{\ in }$\Bbb{R}^{m}$\emph{\ and a tempered distribution }$a$\emph{\
in }$\mathcal{S}_{m}^{\prime }$\emph{\ different from the Dirac delta
distribution }$\delta _{p}$\emph{\ such that the }$p$\emph{-th member of the
family is representable also by} 
\[
v_{p}=\int_{\Bbb{R}^{m}}av. 
\]

\emph{\bigskip }

\emph{Proof.} \textbf{Necessity.} Indeed, if $v_{p}$ fulfills that property
we have that the $^{\mathcal{S}}$linear combination $(a-\delta _{p}).v$ is
zero with a non-zero coefficient distribution, so that the family $v$ is $^{%
\mathcal{S}}$linearly dependent. \textbf{Sufficiency.} Vice versa, let, for
every index-point $p$, the term $v_{p}$ of the family be representable in a
unique way as the superposition $v_{p}=\delta _{p}.v$. Assume $v$ $^{%
\mathcal{S}}$linearly dependent, then there is a coefficient system $a$
different from zero such that $a.v=0$, hence 
\begin{eqnarray*}
v_{p} &=&\int_{\Bbb{R}^{m}}\delta _{p}v-0= \\
&=&\int_{\Bbb{R}^{m}}\delta _{p}v-\int_{\Bbb{R}^{m}}av= \\
&=&\int_{\Bbb{R}^{m}}(\delta _{p}-a)v,
\end{eqnarray*}
since $a$ is a non-zero distribution, the distribution $\delta _{p}-a$ is
different form the Dirac delta $\delta _{p}$, and so the member $v_{p}$ is
representable in another (different) way, against the assumption. $%
\blacksquare $

$\bigskip $

\section{\textbf{Characterizations of }$^{\mathcal{S}}$\textbf{linear
independence}}

\bigskip

By the Dieudonn\'{e}-Schwartz theorem we immediately deduce two
characterizations.

\bigskip

We say an $^{\mathcal{S}}$family $v$ to be \emph{topologically exhaustive}
(with respect to the weak* topology or the strong topology) if its $^{%
\mathcal{S}}$linear hull $^{\mathcal{S}}\mathrm{span}\left( v\right) $\ is $%
\sigma (\mathcal{S}_{n}^{\prime })$-closed (or strongly closed, which is the
same).

\bigskip

\textbf{Theorem.}\emph{\ Let }$v\in \mathcal{S}(\Bbb{R}^{m},\mathcal{S}
_{n}^{\prime })$\emph{\ be a topologically exhaustive family, that is an }$^{%
\mathcal{S}}$\emph{family whose }$^{\mathcal{S}}$\emph{linear hull }$^{%
\mathcal{S}}\mathrm{span}\left( v\right) $\emph{\ is }$\sigma (\mathcal{S}
_{n}^{\prime })$\emph{-closed. Then the following assertions are equivalent}

\begin{itemize}
\item  \emph{1) the family }$v$\emph{\ is }$^{\mathcal{S}}$\emph{linearly
independent;}

\item  \emph{2) the superposition operator }$\int_{\Bbb{R}^{m}}(\cdot ,v)$%
\emph{\ is an injective topological homomorphism for the weak* topologies }$%
\sigma (\mathcal{S}_{m}^{\prime })$\emph{\ and }$\sigma (\mathcal{S}%
_{n}^{\prime })$\emph{;}

\item  \emph{3) the superposition operator }$\int_{\Bbb{R}^{m}}(\cdot ,v)$%
\emph{\ is an injective topological homomorphism for the strong topologies }$%
\beta (\mathcal{S}_{m}^{\prime })$\emph{\ and }$\beta (\mathcal{S}%
_{n}^{\prime })$\emph{;}

\item  \emph{4) the operator }$\widehat{v}$\emph{\ is a surjective
topological homomorphism for the weak topologies }$\sigma (\mathcal{S}_{n})$%
\emph{\ and }$\sigma (\mathcal{S}_{m})$\emph{;}

\item  \emph{5) the operator }$\widehat{v}$\emph{\ is a surjective
topological homomorphism of the topological vector space }$(\mathcal{S}_{n})$%
\emph{\ onto the space }$(\mathcal{S}_{m})$\emph{.}
\end{itemize}

\bigskip

\textbf{Theorem.}\emph{\ Let }$v\in \mathcal{S}(\Bbb{R}^{m},\mathcal{S}
_{n}^{\prime })$\emph{\ be an }$^{\mathcal{S}}$\emph{family. Then the
following assertions are equivalent}

\begin{itemize}
\item  \emph{1) the family }$v$\emph{\ is }$^{\mathcal{S}}$\emph{linearly
independent and the hull }$^{\mathcal{S}}\mathrm{span}\left( v\right) $\emph{%
\ is }$\sigma (\mathcal{S}_{n}^{\prime })$\emph{-closed;}

\item  \emph{2) the superposition operator }$\int_{\Bbb{R}^{m}}(\cdot ,v)$%
\emph{\ is an injective topological homomorphism for the weak* topologies }$%
\sigma (\mathcal{S}_{m}^{\prime })$\emph{\ and }$\sigma (\mathcal{S}%
_{n}^{\prime })$\emph{;}

\item  \emph{3) the operator }$\widehat{v}$\emph{\ is a surjective
topological homomorphism for the weak topologies }$\sigma (\mathcal{S}_{n})$%
\emph{\ and }$\sigma (\mathcal{S}_{m})$\emph{;}

\item  \emph{4) the operator }$\widehat{v}$\emph{\ is a surjective
topological homomorphism from the space }$(\mathcal{S}_{n})$\emph{\ onto the
space }$(\mathcal{S}_{m})$\emph{.}
\end{itemize}

\bigskip

\textbf{Remark (on the coordinate operator).} In the conditions of the above
theorem, if the family $v$ is $^{\mathcal{S}}$linearly independent, we can
consider the algebraic isomorphism from the space $\mathcal{S}_{m}^{\prime }$
onto the $^{\mathcal{S}}$linear hull $^{\mathcal{S}}\mathrm{span}\left(
v\right) $ sending every tempered distribution $a\in \mathcal{S}_{m}^{\prime
}$ into the superposition $a.v$, that is the restriction of the injection 
\[
\int_{\Bbb{R}^{m}}(\cdot ,v):\mathcal{S}_{m}^{\prime }\rightarrow \mathcal{S}%
_{n}^{\prime }:(a,v)\mapsto a.v 
\]
to the pair of sets $(\mathcal{S}_{m}^{\prime },^{\mathcal{S}}\mathrm{span}%
\left( v\right) )$. We shall denote the inverse of this isomorphism by the
symbol $[\cdot |v]$. It is an important consequence of the preceding theorem
that

\begin{itemize}
\item  \textbf{Theorem.} \emph{The operator} 
\[
[\cdot |v]:^{\mathcal{S}}\mathrm{span}\left( v\right) \rightarrow \mathcal{S}%
_{m}^{\prime } 
\]
\emph{is a topological isomorphism, with respect to the topology induced by
the weak* topology }$\sigma (\mathcal{S}_{n}^{\prime })$\emph{\ on the }$^{%
\mathcal{S}}$\emph{linear hull }$^{\mathcal{S}}\mathrm{span}\left( v\right) $%
\emph{\ and to the weak* topology }$\sigma (\mathcal{S}_{m}^{\prime })$\emph{%
, if and only if the }$^{\mathcal{S}}$\emph{linear hull }$^{\mathcal{S}}%
\mathrm{span}\left( v\right) $\emph{\ is }$\sigma (\mathcal{S}_{n}^{\prime
}) $\emph{-closed, that is if the family }$v$\emph{\ is topologically
exhaustive.}
\end{itemize}

$\bigskip $

\section{$^{\mathcal{S}}$\textbf{Bases}}

\bigskip

\textbf{Definition (of }$^{\mathcal{S}}$\textbf{basis).}\emph{\ Let }$v\in 
\mathcal{S}(\Bbb{R}^{m},\mathcal{S}_{n}^{\prime })$\emph{\ be an }$^{%
\mathcal{S}}$\emph{family in }$\mathcal{S}_{n}^{\prime }$ \emph{and let }$V$%
\emph{\ be a subspace of the space }$\mathcal{S}_{n}^{\prime }\emph{.}$\emph{%
\ The family\ }$v$\emph{\ is said an }$^{\mathcal{S}}$\emph{\textbf{basis}
of the subspace }$V$\emph{\ if it is }$^{\mathcal{S}}$\emph{linearly
independent and it }$^{\mathcal{S}}$\emph{generates }$V$\emph{, that is} 
\[
^{\mathcal{S}}\mathrm{span}(v)=V. 
\]
\emph{In other terms, the family\ }$v$\emph{\ is said an }$^{\mathcal{S}}$%
\emph{basis of the subspace }$V$\emph{\ if and only if the superposition
operator of the family }$v$\emph{\ is injective and its range coincides with
the subspace }$V$\emph{.}

\bigskip

The Dirac family $\delta $\ of $\mathcal{S}_{n}^{\prime }$\ is an $^{%
\mathcal{S}}$basis of the whole $\mathcal{S}_{n}^{\prime }$ (indeed its
superposition operator is the identity on $\mathcal{S}_{n}^{\prime }$). We
call the Dirac family $\delta $ \emph{the canonical }$^{\mathcal{S}}$\emph{%
basis of} $\mathcal{S}_{n}^{\prime }$ (because of reasons which will appear
clear soon, and that can be summarized into the relation $u.\delta =u$,
valid for any tempered distribution $u$) or \emph{the Dirac basis of} $%
\mathcal{S}_{n}^{\prime }$.

\bigskip

Moreover, the following complete version of the Fourier expansion-theorem,
allow us to call\ the Fourier families of $\mathcal{S}^{\prime }(\Bbb{R}^{n},%
\Bbb{C})$ by the name of\emph{\ Fourier bases of} $\mathcal{S}^{\prime }(%
\Bbb{R}^{n},\Bbb{C})$.

\bigskip

\textbf{Theorem (geometric form of the Fourier expansion theorem).}\emph{\
In the space of complex tempered distributions }$\mathcal{S}_{n}^{\prime }(%
\Bbb{C})$\emph{\ the Fourier families\ are }$^{\mathcal{S}}$\emph{bases (of
the entire space }$\mathcal{S}_{n}^{\prime }(\Bbb{C})$\emph{).}

\bigskip

\emph{Proof.} Indeed, the superposition operators of the Fourier families in 
$\mathcal{S}_{n}^{\prime }$ are the Fourier transforms upon $\mathcal{S}%
_{n}^{\prime }$ which are bijective. $\blacksquare $

$\bigskip $

\section{\textbf{Algebraic characterizations of }$^{\mathcal{S}}$\textbf{%
bases}}

\bigskip

The following is an elementary but meaningful generalization of the Fourier
expansion theorem.

\bigskip

\textbf{Theorem (characterization of an }$^{\mathcal{S}}$\textbf{basis).}%
\emph{\ Let }$v\in \mathcal{S}(\Bbb{R}^{m},\mathcal{S}_{n}^{\prime })$\emph{%
\ be an }$^{\mathcal{S}}$\emph{family. Then,}

\begin{itemize}
\item  \emph{1) the family }$v$\emph{\ }$^{\mathcal{S}}$\emph{generates the
space }$\mathcal{S}_{n}^{\prime }$\emph{\ if and only if the superposition
operator }$^{t}\left( \widehat{v}\right) $\emph{\ is surjective;}

\item  \emph{2) the family }$v$\emph{\ is }$^{\mathcal{S}}$\emph{linearly
independent if and only if the superposition operator }$^{t}\left( \widehat{v%
}\right) $\emph{\ is injective;}

\item  \emph{3) the family }$v$\emph{\ is an }$^{\mathcal{S}}$\emph{basis of
the space }$\mathcal{S}_{n}^{\prime }$\emph{\ if and only if the
superposition operator }$^{t}\left( \widehat{v}\right) $\emph{\ is bijective.%
}
\end{itemize}

\bigskip

\emph{Proof. }The proof follows immediately from the definitions, however we
see it. First of all the superposition operator $^{t}\left( \widehat{v}%
\right) $ is well defined because $v$ is an $^{\mathcal{S}}$family.
Moreover, it is obvious, by the very definitions, that the family $v$ $^{%
\mathcal{S}}$generates the space $\mathcal{S}_{n}^{\prime }$ if and only if
the superposition operator $^{t}\left( \widehat{v}\right) $ is surjective,
and that $v$ is $^{\mathcal{S}}$linearly independent if and only if the
superposition operator $^{t}\left( \widehat{v}\right) $ is injective.$%
\;\blacksquare $

$\bigskip $

\subsection{\textbf{Example}}

\bigskip

The following point gives us an example of linearly independent family which
is not $^{\mathcal{S}}$linearly independent and also an example of a
linearly independent system of $^{\mathcal{S}}$generators which is not an $^{%
\mathcal{S}}$basis.

\bigskip

\textbf{Example (a system of linearly independent }$^{\mathcal{S}}$\textbf{%
generators that is not an }$^{\mathcal{S}}$\textbf{basis). }Let $v=(\delta
_{x}^{\prime })_{x\in \Bbb{R}}$ be the family in $\mathcal{S}_{1}^{\prime }$
of the first derivatives of the Dirac distributions. The family $v$ is of
class $\mathcal{S}$, in fact 
\[
v(\phi )(x)=v_{x}(\phi )=\delta _{x}^{\prime }(\phi )=-\phi ^{\prime }(x), 
\]
and the the derivative $\phi ^{\prime }$ is an $^{\mathcal{S}}$function.
Consequently, the operator associated with the family $v$ is the derivation
in the test function space $\mathcal{S}_{1}$ up to the sign and, then, the
superposition operator $^{t}\left( \widehat{v}\right) $ is the derivation in
the space $\mathcal{S}_{1}^{\prime }$. This last superposition operator is a
surjective operator (every tempered distribution has a primitive) but it is
not injective (every tempered distribution has many primitives), then the
family $v$ is a system of $^{\mathcal{S}}$generators for the space $\mathcal{%
S}_{1}^{\prime }$, but it is not $^{\mathcal{S}}$linearly independent.
Moreover, note that, however, the family $v$ is linearly independent. In
fact, let $P$ be a finite subset of the real line $\Bbb{R}$, and let, for
every point $p_{0}$ in $P$, $f_{p_{0}}$ be a test function in $\mathcal{S}%
_{1}$ such that 
\[
f_{p_{0}}^{\prime }(p)=\delta _{p_{0},p}, 
\]
for every index $p$ in $P$, here $\delta _{(.,.)}$ is the Kronecker delta
upon the square $P^{2}$. Now, if $a=(a_{p})_{p\in P}$ is a finite family of
scalars such that 
\[
\sum_{p\in P}a_{p}v_{p}=0_{\mathcal{S}_{1}^{\prime }}, 
\]
then 
\[
0=\left( \sum_{p\in P}a_{p}v_{p}\right) (f_{p_{0}})=\sum_{p\in P}a_{p}\delta
_{p_{0}p}=a_{p_{0}}, 
\]
for every index $p_{0}$ in $P$. And so, any finite linear combination of
members of $v$ is zero only with respect to a zero-system of coefficients.

$\bigskip $

\section{\textbf{Totality of }$^{\mathcal{S}}$\textbf{bases}}

\bigskip

We give another way to characterize an $^{\mathcal{S}}$basis.

\bigskip

\textbf{Theorem (characterization of an }$^{\mathcal{S}}$\textbf{basis).}%
\emph{\ Let }$v\in \mathcal{S}(\Bbb{R}^{m},\mathcal{S}_{n}^{\prime })$\emph{%
\ be an }$^{\mathcal{S}}$\emph{family. Then:}

\begin{itemize}
\item[1)]  \emph{the family }$v$\emph{\ is a system of }$^{\mathcal{S}}$%
\emph{generators of the entire space }$\mathcal{S}_{n}^{\prime }$\emph{\ if
and only if the family }$v$\emph{\ is total in the space }$\mathcal{S}_{n}$%
\emph{\ (in the sense that if }$v_{p}(g)=0$\emph{, for every index-point }$p$%
\emph{, then }$g=0$\emph{) and the linear hull of the family }$v$\emph{\ is
weakly* closed;}

\item[2)]  \emph{the family }$v$\emph{\ is }$^{\mathcal{S}}$\emph{linearly
independent if and only it is total in the distribution space }$\mathcal{S}%
_{m}^{\prime }$\emph{, in the sense that if }$a.v=0$\emph{\ then }$a=0$\emph{%
;}

\item[3)]  \emph{a family }$v$\emph{\ is an }$^{\mathcal{S}}$\emph{basis of
the space }$\mathcal{S}_{n}^{\prime }$\emph{\ if and only if the family }$v$%
\emph{\ is total both in the function space }$\mathcal{S}_{n}$\emph{\ and
distribution space }$\mathcal{S}_{m}^{\prime }$\emph{.}
\end{itemize}

\emph{\bigskip }

\emph{Proof.} 1) Indeed, the condition means that the operator generated by
the family $v$ is injective and this (by the Dieudonn\`{e}-Schwartz theorem)
implies the surjectivity of the superposition operator of $v$, since the
image of the superposition operator of $v$ is closed. 2) This is exactly the
definition of linear independence. 3) Follows from the two before. $%
\blacksquare $

$\bigskip $

\section{\textbf{Topological characterizations of }$^{\mathcal{S}}$\textbf{\
bases}}

\bigskip

By the Dieudonn\'{e}-Schwartz theorem we immediately take a characterization.

\bigskip

\textbf{Theorem.}\emph{\ Let }$v\in \mathcal{S}(\Bbb{R}^{m},\mathcal{S}
_{n}^{\prime })$\emph{\ be a family of tempered distributions. Then the
following assertions are equivalent}

\begin{itemize}
\item  \emph{1) the family }$v$\emph{\ is an }$^{\mathcal{S}}$\emph{basis of
the space }$\mathcal{S}_{n}^{\prime }$\emph{;}

\item  \emph{2) the superposition operator }$\int_{\Bbb{R}^{m}}(\cdot ,v)$%
\emph{\ is a topological isomorphism for the weak* topologies }$\sigma (%
\mathcal{S}_{m}^{\prime })$\emph{\ and }$\sigma (\mathcal{S}_{n}^{\prime })$%
\emph{;}

\item  \emph{3) the superposition operator }$\int_{\Bbb{R}^{m}}(\cdot ,v)$%
\emph{\ is a topological isomorphism for the strong topologies }$\beta (%
\mathcal{S}_{m}^{\prime })$\emph{\ and }$\beta (\mathcal{S}_{n}^{\prime })$%
\emph{;}

\item  \emph{4) the operator }$\widehat{v}$\emph{\ is a topological
isomorphism for the weak topologies }$\sigma (\mathcal{S}_{n})$\emph{\ and }$%
\sigma (\mathcal{S}_{m})$\emph{;}

\item  \emph{5) the operator }$\widehat{v}$\emph{\ is a topological
isomorphism of the topological vector space }$(\mathcal{S}_{n})$\emph{\ onto 
}$(\mathcal{S}_{m})$\emph{.}
\end{itemize}

\bigskip

\textbf{David Carf\`{i}}

\emph{Faculty of Economics}

\emph{University of Messina}

\emph{davidcarfi71@yahoo.it}

\bigskip

\end{document}